\documentclass[draft]{article}
\usepackage{amsmath}
\usepackage{amsfonts}
\usepackage{amssymb}
\usepackage[hyphens]{url}

\def \F {\mathbb{F}}
\def \vep {\varepsilon}

\newtheorem{theorem}{Theorem}
\newtheorem{lemma}{Lemma}

\newtheorem{proposition}{Proposition}

\begin{document}

\title{Expansion Complexity and Linear Complexity 
 of Sequences over Finite Fields} 

\author{ 
L\'aszl\'o M\'erai, Harald Niederreiter, Arne Winterhof\\
 {\small Johann Radon Institute for Computational and Applied Mathematics}\\
 {\small Austrian Academy of Sciences}\\
 {\small Altenbergerstr.\ 69,
 4040 Linz, Austria}\\
{\small  e-mail: \{\texttt{laszlo.merai,harald.niederreiter,arne.winterhof}\}\texttt{@oeaw.ac.at}}
 } 
 
\maketitle
 
 \begin{abstract}
The linear complexity is a measure for the unpredictability of a sequence over a finite field and thus for its suitability in cryptography. 
 In 2012, Diem introduced a new figure of merit for cryptographic sequences called expansion complexity. 
 We study the relationship between linear complexity and expansion complexity. In particular, we show that 
 for purely periodic sequences both figures of merit provide essentially the same quality test for a sufficiently long part of the sequence. 
 However, if we study shorter parts of the period or nonperiodic sequences, then
 we can show, roughly speaking, that the expansion complexity provides a stronger test. We demonstrate this by analyzing a sequence of binomial coefficients modulo $p$.
 Finally, we establish a probabilistic result on the behavior of the expansion complexity of random sequences over a finite field.
\end{abstract}

 
 
\let\thefootnote\relax\footnote{The final publication is available at Springer via \url{http://dx.doi.org/10.1007/s12095-016-0189-2}}

\section{Introduction}

For a sequence ${\cal S}=(s_i)_{i=0}^\infty$ over the finite field $\F_q$ of $q$ elements and a positive integer $N$, the {\em $N$th linear complexity} $L_N=L_N({\cal S})$
is the length of a shortest linear recurrence
\begin{equation}\label{rec}
 s_{i+L_N}+\sum_{\ell=0}^{L_N-1} c_\ell s_{i+\ell}=0, \quad 0\le i\le N-L_N-1,
\end{equation}
with coefficients $c_\ell\in \F_q$, which is satisfied by the first $N$ terms of the sequence. We use the convention $L_N=0$ if $s_0=s_1=\ldots=s_{N-1}=0$
and $L_N=N$ if $s_0=s_1=\ldots=s_{N-2}=0\ne s_{N-1}$. 
The {\em linear complexity} $L=L({\cal S})$ is
$$L({\cal S})=\sup_{N\ge 1} L_N({\cal S}).$$
Note that $L$ is finite if and only if ${\cal S}$ is ultimately periodic. If $T$ and $t$ denote the period and preperiod of ${\cal S}$, respectively, we obviously have
$$L\le T+t.$$

The ($N$th) linear complexity is a measure for the unpredictability of a sequence and thus its suitability in cryptography. A sequence with small $L_N$ for a sufficiently large $N$
is disastrous for cryptographic applications. However, the converse is not true. There are highly predictable sequences with large $L_N$, including the example 
$s_0=\ldots=s_{N-2}=0\ne s_{N-1}$. Hence, for testing the suitability of a sequence in cryptography we also have to study finer figures of merit. A recent survey on the
linear complexity is given in~\cite{MW13}.

Diem \cite{di12} introduced the expansion complexity of the sequence ${\cal S}$ as follows. We define the {\em generating function} $G(x)$ of ${\cal S}$ by
$$G(x)=\sum_{i=0}^\infty s_i x^i,$$
viewed as a formal power series over $\F_q$. Note the change by the factor $x$ compared to the definition in~\cite{di12}.
For a positive integer $N$, the {\em $N$th expansion complexity} $E_N=E_N({\cal S})$ is $E_N=0$ if $s_0=\ldots=s_{N-1}=0$ and otherwise the least total degree
of a nonzero polynomial $h(x,y)\in \F_q[x,y]$ with
\begin{equation} \label{eqcon}
h(x,G(x))\equiv 0 \bmod x^N.
\end{equation}
Note that $E_N$ depends only on the first $N$ terms of ${\cal S}$.

We prove upper and lower bounds on $E_N$ in terms of $L_N$ and the smallest number $t_N$ with $c_{t_N}\ne 0$ in $(\ref{rec})$.  
In particular, we show that 
 for purely periodic sequences both figures of merit provide essentially the same quality test for the whole sequences. However, if we study only parts of the period or nonperiodic sequences,
 we can show, roughly speaking, that the expansion complexity provides a stronger test. We demonstrate this by analyzing linear complexity and expansion complexity of the sequence
 ${\cal A}=(a_i)_{i=0}^\infty$ of binomial coefficients $a_i={i+k\choose k}$ modulo a prime $p$ for some $1\le k\le p-1$.

First we study ultimately periodic sequences in Section~\ref{periodseq}. Then we analyze the linear complexity and expansion complexity of the sequence ${\cal A}$ in Section~\ref{bin}.
The aperiodic case is studied in Section~\ref{aper}. A probabilistic result on the behavior of the expansion complexity of random sequences over a finite field is shown in
Section~\ref{sepr}. 
 
\section{Ultimately periodic sequences}
\label{periodseq}
Now let ${\cal S}=(s_i)_{i=0}^\infty$ be an ultimately periodic sequence over $\F_q$ with preperiod~$t$ and period $T$, that is, $s_{i+t+T}=s_{i+t}$ for $i=0,1,\ldots$ . Let $L$ be its 
linear complexity and recall that $L\le T+t$.  
Then its generating function is a rational function
\begin{equation}\label{gen}
  G(x)=\frac{f(x)}{g(x)}
\end{equation}
with polynomials $f(x),g(x)\in \F_q[x]$ with $\deg(f)<L$, $\deg(g)=L-t$,  and $\gcd(f(x),g(x))=\gcd(g(x),x)=1$,
see \cite[Theorem 8.40]{lini83}. Note that such a sequence satisfies a linear recurrence of the form
$$\sum_{\ell=t}^Lc_\ell s_{i+\ell}=0,\quad i\ge 0,$$
with $c_L=1$ and $c_t\ne 0$. Then we have
\begin{equation}\label{gx}
g(x)=1+c_{L-1}x+\cdots+c_tx^{L-t}.
\end{equation}

 \begin{lemma}\label{noN}
 Let $G(x)$ in~\eqref{gen} be not identically zero
 and let $h(x,y)\in \F_q[x,y]$ be a nonzero polynomial of local degree $d$ in~$y$.
 Put $H(x)=g(x)^dh(x,G(x))$. If $H(x)$ is the zero polynomial, then the total degree of $h(x,y)$ satisfies
 $$\deg(h)\ge L-t+1.$$
\end{lemma}

Proof. 
We write 
$$h(x,y)=\sum_{i=0}^d h_i(x)y^i\in \F_q[x,y]$$  
with $h_d(x)\ne 0$. Then $H(x)=0$ implies
\begin{equation}\label{fg}
  \sum_{i=0}^d h_i(x)f(x)^i g(x)^{d-i}=0 
\end{equation}
and $d\ge 1$, where we used $(\ref{gen})$. Note that $g(x)\ne 0$ by $(\ref{gx})$. 
Hence, $h_d(x)$ is divisible by $g(x)$ and thus of degree at least $\deg(g)=L-t$. 
Finally, we have $\deg(h)\ge \deg(h_d)+d\ge L-t+1$.~\hfill $\Box$

\begin{theorem}
\label{periodbound}
 Let ${\cal S}$ be an ultimately periodic sequence over $\F_q$ with preperiod $t$, linear complexity $L$, and generating function $G(x)\ne 0$.
 Then we have
 $$ E_N({\cal S})\ge \left\{\begin{array}{ll} L-t+1 & \mbox{for } N>(L-t)(L-\min\{1,t-1\}),\\
                       \lceil N/(L-\min\{1,t-1\})\rceil & \mbox{otherwise,}
                       \end{array}\right.$$
 and
 $$E_N({\cal S})\le L+\max\{-1,-t+1\}.$$
\end{theorem}

Proof. Since otherwise the lower bound is trivial, we may assume $\deg(h)< N/(L-\min\{1,t-1\})$. Then $\deg(H)\le \deg(h)(L-\min\{1,t-1\})<N$ using $(\ref{fg})$ and  
$$h(x,G(x))\equiv 0 \bmod x^N$$
is equivalent to $H(x)=0$. Now the lower bound follows by Lemma~\ref{noN}. 

Choosing the polynomial
$$h(x,y)=g(x)y-f(x)$$
of degree $\deg(h)=\max\{\deg(f),\deg(g)+1\}\le \max\{L-1,L-t+1\}$,
we get the upper bound.~\hfill $\Box$ \\

Remark.  For $t\le 2$ and $N>(L-t)(L-t+1)$ we have equality: 
$$E_N({\cal S})=L-t+1.$$

\section{A sequence of binomial coefficients}
\label{bin}

For a prime $p$ and some integer $k$ with $1\le k\le p-1$, we study the $p$-periodic sequence ${\cal A}=(a_i)_{i=0}^\infty$ of binomial coefficients
\begin{equation} \label{bindef} 
a_i={i+k \choose k}\bmod p,\quad i=0,1,\ldots .
\end{equation}
First we will show that ${\cal A}$ has an optimal $N$th linear complexity for $1\le N\le 2\min\{k+1,p-k\}$ which suggests an optimal value of $k=(p-1)/2$.
However, since the last $k$ sequence elements $a_{p-k},a_{p-k+1},\ldots,a_{p-1}$ in the first period vanish, the sequence becomes more predictable with increasing $k$.  

It turns out that the $p$th expansion complexity is $E_p({\cal A})=\min\{k+2,\lceil p/(k+2)\rceil\}$ which suggests an optimal value of $k\approx p^{1/2}$, where only the 
first $p-k$ sequence elements should be used in practice. 

\subsection{Linear complexity}

\begin{proposition}
\label{lincompA}
 We have 
 $$L({\cal A})=k+1$$
 and
 $$L_N({\cal A})\ge \min\{k+1,\lceil N/2\rceil,p-k\}.$$
\end{proposition}

Proof. 
Since ${i+k \choose k}=\prod_{j=1}^k \frac{i+j}{j}$ is a polynomial of degree $k$ in $i$, we can apply 
the following well-known result, see
\cite[Theorem~8]{bletpa96} or \cite[Theorem~1]{mewi02}, to get the value of the linear complexity:
let $f$ be a polynomial of degree $d<p$ over $\F_p$ and ${\cal S}=(s_i)_{i=0}^\infty$ be the $p$-periodic sequence defined by
$s_i=f(i)$ for $i=0,1,\ldots$; 
then $L({\cal S})=d+1$. Furthermore, $L_N({\cal S})\ge \min\{d+1,N-d\}$ by \cite[Theorem~3]{mewi03}, which implies
$$L_N({\cal A})\ge \min\{k+1,N-k\}.$$

Put $L=L_N({\cal A})$. Since otherwise the second result is trivial, we may assume
$$L\le \min\{k,p-k-1\}\quad \mbox{and}\quad N\le \min\{2k,p-1\}.$$
Assume there is a linear recurrence of length $L$ satisfied by the first $N$ terms of~${\cal A}$, that is,
$$\sum_{\ell=0}^{L} c_\ell a_{i+\ell}=0,\quad 0\le i\le N-L-1,$$
where $c_L=-1$.
Note that 
$$a_{i+\ell}={i+\ell +k \choose k}= a_i \prod_{j=1}^\ell \frac{i+k+j}{i+j}.$$
With $f_\ell(x)=\prod_{j=1}^\ell (x+k+j)$ and $g_\ell(x)=\prod_{j=1}^\ell (x+j)$,
we get
$$\sum_{\ell=0}^L c_\ell \frac{f_\ell(i)}{g_\ell(i)}=0,\quad 0\le i\le \min\{N-L,p-k\}-1,$$
since $a_i\ne 0$ for $0\le i\le p-k-1$. Multiplying with $g_L(i)$, we get
$$\sum_{\ell=0}^L c_\ell f_\ell(i)\prod_{j=\ell+1}^L(i+j)=0,\quad 0\le i\le \min\{N-L,p-k\}-1.$$
We have constructed a polynomial of degree at most $L$ with at least $\min\{N-L,p-k\}$ zeros.
Evaluating the left hand side at $i=p-L\ge p-k$, we get the value $c_Lf_L(p-L)\ne 0$. Hence by Lagrange's theorem we obtain
$$L\ge \min\{N-L,p-k\}.$$
If $L\ge N-p+k$, we get $L\ge \max\{N/2,N-p+k\}=N/2$ since
$N-p+k\le L<p-k$ implies $N<2(p-k)$.
If $L<N-p+k$, we obtain $L\ge p-k$.
~\hfill $\Box$

\subsection{Expansion complexity}

\begin{lemma}
\label{genbin}
 The generating function $G(x)$ of ${\cal A}$ is
 $$G(x)=\frac{1}{(1-x)^{k+1}}.$$
\end{lemma}

Proof. 
First verify that 
$${p-1-k\choose i}(-1)^i\equiv \prod_{j=1}^i \frac{k+j}{j}\equiv {i+k \choose i}\equiv {i+k\choose k}\bmod p.$$
Then we get 
\begin{eqnarray*}
(1-x)^pG(x)&=&(1-x^p)G(x)=\sum_{i=0}^{p-1-k} {i+k \choose k}x^i\\
&=&\sum_{i=0}^{p-1-k} {p-1-k \choose i}(-x)^i=(1-x)^{p-1-k} 
\end{eqnarray*}
and the result follows.~\hfill $\Box$

\begin{theorem}
 Let ${\cal A}=(a_i)_{i=0}^\infty$ be the sequence of binomial coefficients modulo $p$ defined by $(\ref{bindef})$ and $E_p({\cal A})$ its $p$th expansion complexity.
 
 For $(k+1)(k+2)<p$ we have
 \[
  E_p({\cal A})=k+2
 \]
and for $(k+1)(k+2)\geq p$ 
\[
 \left\lceil \frac{p}{k+2} \right\rceil \leq E_p({\cal A}) \leq \max\left\{\left\lceil \frac{p}{k+2}\right\rceil,  (k+1)\left\{\frac{p}{k+1}\right\}\right\},
\]
where $\{x\}$ is the fractional part of $x$, that is,  $\{x\}=x-\lfloor x \rfloor$.
\end{theorem}

Proof. 
By Proposition~\ref{lincompA} we have $L=L({\cal A})=k+1$. If $(k+1)(k+2)<p$ we get by Theorem~\ref{periodbound} (with $t=0$ since ${\cal A}$ is purely periodic) the first result.

If $(k+1)(k+2) \geq p$ we have by Theorem~\ref{periodbound}
\[
  E_p({\cal A}) \geq \left\lceil \frac{p}{k+2} \right\rceil .
\]
We put 
$$d=\min\left\{\left\lfloor \frac{p}{k+1}\right\rfloor, \left\lceil \frac{p}{k+2}\right\rceil\right\}$$
and take 
$$h(x,y)=y^d-(1-x)^{p-d (k+1)}\in\mathbb{F}_p[x,y].$$
Here we used $d\leq p/(k+1)$ since otherwise $h(x,y)$ is not a polynomial.
By Lemma~\ref{genbin}
we have $G(x)=\frac{1}{(1-x)^{k+1}}$ and thus
\begin{eqnarray*} 
h(x,G(x))&=&\frac{1}{(1-x)^{d(k+1)}}-(1-x)^{p-d(k+1)}
=\frac{1-(1-x)^p}{(1-x)^{d(k+1)}}\\
&=&\frac{x^p}{(1-x)^{d(k+1)}}\equiv 0 \bmod x^p
\end{eqnarray*}
since $\gcd((1-x),x)=1$.
Hence, 
$$
E_p({\cal A})\le \deg(h)=\max\{d,p-d(k+1)\}=
\begin{cases}
d & \text{if }  d= \left\lceil \frac{p}{k+2}\right\rceil,\\
p-d(k+1) & \text{otherwise,}
\end{cases}
$$
and the result follows.~\hfill $\Box$

\section{The aperiodic case}
\label{aper}

\subsection{Growth of $E_N({\cal S})$ and $L_N({\cal S})$}

First we describe the possible growth of the nondecreasing function\\ $N\mapsto~E_N({\cal S})$.

\begin{proposition}\label{engrowth}
 We have $E_N({\cal S})\le E_{N+1}({\cal S})\le E_N({\cal S})+1$.
\end{proposition}

Proof. If $h(x,G(x))\equiv 0 \bmod x^N$, then $xh(x,G(x))\equiv 0\bmod x^{N+1}$.~\hfill $\Box$\\

For comparison, we state the corresponding result on the possible growth of the nondecreasing function $N \mapsto L_N({\cal S})$, which is called the {\em linear complexity profile}
of ${\cal S}$. For a proof see \cite[Theorem 6.7.4]{ju93}, \cite{ma69}, or \cite[Chapter~4]{ru86}.

\begin{lemma}\label{lngrowth}
 If $L_N({\cal S}) > N/2$, then $L_{N+1}({\cal S}) =
L_N({\cal S})$. 
If $L_N({\cal S})\le N/2$, then $L_{N+1}({\cal S})\in \{L_N({\cal S}),N + 1 - L_N({\cal S})\}$.
\end{lemma}

\subsection{Bounds}

\begin{theorem}
\label{enln}
Let ${\cal S}$ be a sequence over $\F_q$ with generating function $G(x)$.
For $N\ge 2$ let $G(x)$ satisfy 
$$G(x)\not\equiv 0\bmod x^N$$
and let 
$$
\sum_{\ell=t_N}^{L_N} c_\ell s_{i+\ell}=0,\quad 0\le i\le N-L_N-1, 
$$
be a shortest linear recurrence for the first $N$ terms of ${\cal S}$, where $c_{L_N}=1$ and $c_{t_N}\ne 0$.
 Then 
 $$E_N({\cal S})\ge \left\{\begin{array}{ll} L_N-t_N+1 & \mbox{for } N>(L_N-t_N)(L_N-\min\{1,t_N-1\}),\\
                       \left\lceil \frac{N}{L_N-\min\{1,t_N-1\}}\right\rceil & \mbox{otherwise,}
                       \end{array}\right.$$
and                       
$$E_N({\cal S})\le \min\{L_N({\cal S})+\max\{-1,-t_N+1\},N-L_N({\cal S})+2\}.$$ 
\end{theorem}

Proof.
Let ${\cal U}=(u_i)_{i=0}^\infty$ be the ultimately periodic sequence with preperiod $t=t_N$ defined by 
$$u_i=s_i \mbox{ for }i=0,1,\ldots,N-1$$ 
and 
$$u_{i+L_N}=-\sum_{\ell=t}^{L_N-1}c_\ell s_{i+\ell} \mbox{ for }i=N-L_N,N-L_N+1,\ldots .$$
Then we have $E_N({\cal S})=E_N({\cal U})$ and $L_N({\cal S})=L_N({\cal U})=L({\cal U})$.  
By Theorem~\ref{periodbound}, 
it remains to show that $E_N\le N-L_N+2$ if $L_N>(N+1)/2$.
In particular, we have already proved that
\begin{equation}\label{firstbound}
 E_N({\cal S})\le L_N({\cal S})+1.
\end{equation}

If $L_N({\cal S})> (N+1)/2$, then we have $L_N({\cal S})=L_{N-1}({\cal S})=\ldots=L_{N-k}({\cal S})\ne L_{N-k-1}({\cal S})$
for some $0\le k<(N-1)/2$ since $(N+1)/2<L_N({\cal S})=L_{N-k}({\cal S})\le N-k$.
By Lemma~\ref{lngrowth} we get
$L_N({\cal S})=L_{N-k}({\cal S})=N-k-L_{N-k-1}({\cal S})\le N-k-E_{N-k-1}({\cal S})+1\le N-E_N({\cal S})+2$
by $(\ref{firstbound})$ and Proposition~\ref{engrowth},
and the remaining bound follows. 
~\hfill $\Box$\\

Remarks.
\begin{itemize}
 \item For $N\ge 2$ we have 
 $$E_N({\cal S})\le \min\left\{\left\lfloor \frac{N+3}{2}\right\rfloor,N-1\right\},$$
 where $E_N\le N-1$ can be obtained 
 by choosing $h(x,y)=y-\sum\limits_{i=0}^{N-1}s_iX^i$ in~\eqref{eqcon}.
 \item The Berlekamp-Massey algorithm does not only compute the whole linear complexity profile $L_N({\cal S})$ for $N=1,2,\ldots$, 
 but also shortest linear recurrences satisfied by the first $N$ terms, from which we can get $t_N$ as well; see for example \cite{be68,ju93,ma69}. 
 \item We may modify Diem's definition by adding the condition that $h(x,y)$ is irreducible over $\F_q$.
Without this modification $E_N({\cal S})$ may depend only on the first $N_0$ terms of ${\cal S}$ for some $N_0<N$. 
More precisely, assume that all $h(x,y)\ne 0$ of minimal degree satisfying~\eqref{eqcon}
are of the form
$h(x,y)=h_1(x,y)h_2(x,y)$ with nonconstant polynomials $h_1(x,y)$ and $h_2(x,y)$ over~$\F_q$.  
Then $h_1(x,G(x))\equiv 0\bmod x^{N_1}$ and
$h_2(x,G(x))\equiv 0 \bmod x^{N_2}$ for some $1\le N_1,N_2<N$ with $N=N_1+N_2$, and so
$E_N({\cal S})$ depends only on the  first $N_0=\max\{N_1,N_2\}$ terms of ${\cal S}$.
However, using only irreducible polynomials would cause serious modifications in the algorithm suggested in \cite[Section~5]{di12}.
\item We have $E_{N_1+N_2}\le E_{N_1}+E_{N_2}$ if $G(x)\not\equiv 0 \bmod x^{\min\{N_1,N_2\}}$.
Indeed, if $h_1(x,G(x))\equiv 0\bmod x^{N_1}$ and $h_2(x,G(x))\equiv 0 \bmod x^{N_2}$ with nontrivial polynomials $h_1(x,y)$ and $h_2(x,y)$, then $h(x,G(x))\equiv 0\bmod x^{N_1+N_2}$. 
\item Let $p$ be the characteristic of $\F_q$. For $N\ge 2$ let $k$ be the nonnegative integer with $p^k\le N-1<p^{k+1}$.
Then we have $E_N\le \lfloor (N-1)/p^k\rfloor p^k$ taking $h(x,y)=y^{p^k}-\sum_{i=0}^{\lfloor (N-1)/p^k\rfloor}s_ix^{ip^k}$, which improves
$E_N\le (N+3)/2$ in some cases.
\end{itemize}

\section{A probabilistic result}
\label{sepr}

Let $\mu_q$ be the uniform probability measure on $\F_q$ which assigns the measure $1/q$ to each element of $\F_q$. Let $\F_q^{\infty}$ be the sequence space over $\F_q$
and let $\mu_q^{\infty}$ be the complete product probability measure on $\F_q^{\infty}$ induced by $\mu_q$. We say that a property of sequences ${\cal S} \in \F_q^{\infty}$
holds $\mu_q^{\infty}$-\emph{almost everywhere} if it holds for a set of sequences ${\cal S}$ of $\mu_q^{\infty}$-measure $1$. We may view such a property as a typical
property of a random sequence over $\F_q$.

\begin{theorem}
\label{thpr}
We have
$$
\liminf_{N \to \infty} \, \frac{E_N({\cal S})}{N^{1/2}} \ge 1 
\qquad \mu_q^{\infty}\mbox{-almost everywhere}.
$$
\end{theorem}

Proof. First we fix $\vep$ with $0 < \vep < 1$ and we put
\begin{equation} \label{eqbn}
b_N=\left\lfloor (1-\vep)^{1/2} N^{1/2} \right\rfloor \qquad \mbox{for } N=1,2,\ldots .
\end{equation}
Then $b_N \ge 1$ for all sufficiently large $N$. For such $N$ put
$$
A_N=\{{\cal S} \in \F_q^{\infty}: E_N({\cal S}) \le b_N\}.
$$
Since $E_N({\cal S})$ depends only on the first $N$ terms of ${\cal S}$, the measure $\mu_q^{\infty}(A_N)$ is given by
\begin{equation} \label{equq}
\mu_q^{\infty}(A_N) =q^{-N} \cdot \# \{{\cal S} \in \F_q^N: E_N({\cal S}) \le b_N\}.
\end{equation} 
According to \cite[Proposition 7]{di12}, ${\cal S}$ is uniquely determined by its
first
$b_N^2$ terms. It follows therefore that
$$
\# \{{\cal S} \in \F_q^N: E_N({\cal S}) \le b_N\} \le q^{b_N^2}. 
$$
It follows thus from~\eqref{eqbn} and~\eqref{equq} that $\mu_q^{\infty}(A_N) \le q^{-\vep N}$ for all sufficiently large $N$. Therefore
$\sum_{N=1}^{\infty} \mu_q^{\infty}(A_N) < \infty$. Then the Borel-Cantelli lemma (see \cite[Lemma~3.14]{Br} and \cite[p.~228]{Lo}) shows that the set of all ${\cal S} \in
\F_q^{\infty}$ for which ${\cal S} \in A_N$ for infinitely many $N$ has $\mu_q^{\infty}$-measure $0$. In other words, $\mu_q^{\infty}$-almost everywhere we have
${\cal S} \in A_N$ for at most finitely many $N$. It follows then from the definition of $A_N$ that $\mu_q^{\infty}$-almost everywhere we have
$$
E_N({\cal S}) > b_N > 
(1-\vep)^{1/2} N^{1/2} -1
$$
for all sufficiently large $N$. Therefore $\mu_q^{\infty}$-almost everywhere,
$$
\liminf_{N \to \infty} \, \frac{E_N({\cal S})}{N^{1/2}} \ge 
(1-\vep)^{1/2}.
$$
By applying this for $\vep =1/r$ with $r=1,2,\ldots$ and noting that the intersection of countably many sets of $\mu_q^{\infty}$-measure $1$ has again $\mu_q^{\infty}$-measure $1$,
we obtain the result of the theorem.~\hfill $\Box$\\

Theorem~\ref{thpr} shows that, for random sequences ${\cal S}$ over $\F_q$, the expansion complexity $E_N({\cal S})$ grows at least at the rate $N^{1/2}$ as $N \to \infty$.
It may be conjectured that this is the exact order of magnitude of $E_N({\cal S})$ for random sequences~ ${\cal S}$ over~$\F_q$.

\section*{Acknowledgements}
The authors wish to thank Claus Diem for a hint which led to an improvement of the constant in Theorem~\ref{thpr}.

The first and the third author are partially supported by the Austrian Science Fund FWF Project F5511-N26 
which is part of the Special Research Program "Quasi-Monte Carlo Methods: Theory and Applications".

\end{document}